\documentclass[10pt]{article}
\title{\textit {Derived invariance of Hochschild-Mitchell (co)homology 
and one-point extensions}}

\author{Estanislao Herscovich and Andrea Solotar
\thanks{This work has been supported by the projects PICT 08280 
(ANPCyT), UBACYTX169 and PIP-CONICET 5099.
The first author is a CONICET fellow.
The second author is a research member of CONICET (Argentina) and a 
Regular
Associate of ICTP Associate Scheme.}}

\date{}

\input {xy}
\usepackage{amsmath,amsthm,amsfonts,amssymb}
\usepackage[matrix,arrow]{xy}

\usepackage{graphicx}
\newtheorem{teo}{Theorem}[section]
\newtheorem{coro}[teo]{Corollary}
\newtheorem{lema}[teo]{Lemma}
\newtheorem{prop}[teo]{Proposition}
\newtheorem{defi}[teo]{Definition}
\newtheorem{obs}[teo]{Remark}
\newtheorem{rem}[teo]{Remark}
\newtheorem{ejem}[teo]{Example}

\usepackage{fullpage}

\numberwithin{equation}{section}     

\def\id{{\mathrm{id}}}
\def\cl#1{{\langle #1 \rangle}}
\def\QED{\hfil$\Box$}
\def\qed{\begin{flushright} \QED \end{flushright}}

\newcommand\ZZ{{\mathbb{Z}}}
\newcommand\NN{{\mathbb{N}}}

\def\C{{\mathcal C}}
\def\D{{\mathcal D}}
\def\E{{\mathcal E}}
\def\F{{\mathcal F}}
\def\H{{\mathcal H}}

\def\place{{-}}

\begin{document}
\maketitle

\sf

\begin{abstract}
In this article we prove derived invariance of Hochschild-Mitchell 
homology 
and cohomology 
and we extend to $k$-linear categories a result by Barot and Lenzing 
concerning derived equivalences and one-point extensions. 

We also prove the existence of a long exact sequence \`a la Happel and 
we give a generalization of this result which provides an alternative 
approach. 
\end{abstract}
                                                                                
                                                                                

\normalsize                                                                                

\section{Introduction}

It is known that linear categories over a field $k$ are a 
generalization of finite dimensional $k$-algebras: 
given a finite dimensional unitary $k$-algebra $A$ and a complete 
system $E =\{ e_{1}, \dots, e_{n}\}$ of orthogonal 
idempotents of $A$, the category $\C_{A}$ with objects indexed by $E$ 
and morphisms given by 
$\mathrm{Hom}_{\C}(e_{i},e_{j}) = e_{j} A e_{i}$ may be associated to 
$A$. 
Different complete sets of orthogonal idempotents of $A$ give different 
categories, but all of them are 
Morita equivalent \cite{CS1}. 
Conversely given a $k$-linear category with a finite set of objects 
$\C_{0} =\{x_{1}, \dots, x_{n}\}$, 
$a(\C) = \oplus_{i,j = 1}^{n} \mathrm{Hom}_{\C} (x_{i} , x_{j})$ is a 
$k$-algebra with unit 
$\sum_{i=1}^{n} \id_{x_{i}}$. 
The categorical point of view gives in our opinion a very clear 
insight. 

Hochschild cohomology is a very powerful tool in the study of finite dimensional algebras. 
Its counterpart when working with $k$-linear categories is Hochschild-Mitchell 
cohomology \cite{Mit2}. 
Some computations of Hochschild-Mitchell cohomology groups have been achieved recently by 
de la Pe\~na and Clotilde Garc\'{\i}a \cite{dP-G1} and by ourselves \cite{HS1}. 

In this article we study one-point extensions of linear $k$-categories, 
obtaining two main results. 
The first one concerns derived invariance of Hochschild-Mitchell 
cohomology, and the second one is the existence 
of a cohomological long exact sequence relating the cohomology of the 
category itself and the cohomology of its one-point extension. 
Both of them should be useful for computations, as it is the case for algebras.

More precisely, let $A$ be a finite dimensional $k$-algebra and $M$ a 
right $A$-module. 
It has been proved by Barot and Lenzing \cite{BL1} that if $A$ is 
derived equivalent to another finite 
dimensional $k$-algebra $B$ and the equivalence maps $M$ into a right 
$B$-module $N$, then the one point extensions 
$A[M]$ and $B[N]$ are such that there exists a triangulated equivalence 
$\Phi : D^{b}(\mathrm{Mod}_{A[M]}) \rightarrow 
D^{b}(\mathrm{Mod}_{B[N]})$ and $\Phi$ restricts to a triangulated 
equivalence $\phi : D^{b}(\mathrm{Mod}_{A}) \rightarrow 
D^{b}(\mathrm{Mod}_{B})$. 
The motivation of this article was to prove that this result holds for 
small $k$-linear categories instead of finite dimensional $k$-algebras. 
This is achieved in Theorem \ref{thm:derone}. 

In the way to prove this theorem, we give in Theorem \ref{thm:2} an 
alternative description of Morita 
equivalences between $k$-linear categories (cf. \cite{CS1}) and a 
characterization of derived equivalences in this context (Theorem 
\ref{thm:Kel9.2}). 
As a consequence we prove in Theorem \ref{thm:hoch} that 
Hochschild-Mitchell homology and cohomology are derived invariant. 

We also prove in section $4$ that the Hochschild-Mitchell cohomology of a one point 
extension is related to the 
Hochschild-Mitchell cohomology of the category by a long exact sequence 
\textit{\`a la Happel} \cite{Hap1}. 
Actually, we prove this fact in two different ways. 
Firstly, we provide a direct proof and secondly we reobtain the result 
as an example of a much more general situation (cf. Thm. 
\ref{thm:longseq}). 
The analogue for finite dimensional algebras is proved in \cite{Cib1} 
and \cite{MP1}. 
Our proof is related to Cibils' article, but it is fact simpler, even for 
the case of algebras. 

We thank M. Su\'arez-\'Alvarez for useful discussions. 

\section{Morita and derived equivalences}

In the first part of this section we shall give a description of 
equivalences between the module categories 
of two $k$-linear categories $\C$ and $\D$ that will lead to a 
characterization of Morita equivalences which 
is in fact very close to the algebraic case. 

We begin by recalling the definition of a module over a linear category 
$\C$. 
For further references, see  \cite{Mit2}, \cite{CM1} and \cite{CR1}.    

Let us consider a field $k$ and a small category $\C$. 
 
$\C$ is a \textbf{$k$-linear category} if the set of morphisms between 
two arbitrary objects of $\C$ is a $k$-vector space and 
composition of morphisms is $k$-bilinear. 
From now on, $\C$ will be a $k$-linear category with set of objects  
$\C_{0}$ and given objects $x, y$ we shall denote ${}_{y}\C_{x}$ the 
$k$-vector space of morphisms from $x$ to $y$ in $\C$. 
Given $x,y,z$ in $\C_{0}$, the composition is a $k$-linear map 
\[     \circ_{z,y,x} : {}_{z}\C_{y} \otimes {}_{y}\C_{x} \to 
{}_{z}\C_{x}.     \]
We shall denote ${}_{z}f_{y}.{}_{y}g_{x}$, 
${}_{z}f_{y}\cdot{}_{y}g_{x}$, ${}_{z}f_{y} \circ {}_{y}g_{x}$, or $f g$ if subscripts are clear, 
the image of 
${}_{z}f_{y} \otimes {}_{y}g_{x}$ under this map.

The simplest example of $k$-linear category is to look at a $k$-algebra 
$A$ as a category with only one object and the set of morphisms equal 
to $A$.  

\begin{defi}
Given two $k$-linear categories $\C$ and $\D$ \textbf{the (external) 
tensor product category}, which we denote $\C \boxtimes_{k} \D$, is the 
category with 
set of objects $\C_{0} \times \D_{0}$ and given $c, c' \in \C_{0}$ and 
$d,d' \in \D_{0}$ 
\[     {}_{(c',d)}(\C \boxtimes_{k} \D)_{(c,d)} = {}_{c}\C_{c} \otimes 
{}_{d'}\D_{d}.       \]

The functor $\C \boxtimes_{k} \place$ is the left adjoint functor to 
$\mathrm{Func}(\C, \place)$ (see \cite{Mit2}, section 2, p. 13). 
We will omit the subindex $k$ in the external tensor product. 
We will call the category $\C \boxtimes \C^{op}$ the \textbf{enveloping 
category of $\C$} and denote it $\C^{e}$. 
\end{defi}

\begin{defi}
A \textbf{left $\C$-module} $M$ is a covariant $k$-linear functor from 
the category $\C$ to the category of vector spaces over $k$. 
Equivalently, a left $\C$-module $M$ is a collection of $k$-vector 
spaces $\{ {}_{x}M\}_{x \in \C_{0}}$ provided with a left action 
\[     {}_{y}\C_{x} \otimes {}_{x}M \rightarrow {}_{y}M,     \]
where the image of ${}_{y}f_{x} \otimes {}_{x}m$ is denoted by 
${}_{y}f_{x}.{}_{x}m$ or $f m$, satisfying the usual axioms 
\begin{gather*}
   {}_{z}f_{y}.({}_{y}g_{x}.{}_{x}m) = 
({}_{z}f_{y}.{}_{y}g_{x}).{}_{x}m,
   \\
   {}_{x}1_{x}.{}_{x}m = {}_{x}m.
\end{gather*}

Right $\C$-modules are defined in an analogous way. 
Also, a \textbf{$\C$-bimodule} is just a $\C^{e}$-module. 

We shall denote ${}_{\C}\mathrm{Mod}$ and $\mathrm{Mod}_{\C}$ the 
categories of left $\C$-modules and right 
$\C$-modules, respectively. 
\end{defi}

The obvious example of $\C$-bimodule is given by the category itself, 
i.e. ${}_{y}\C_{x}$ for every $x,y \in \C_{0}$. 
We will denote this bimodule by $\C$. 

In a similar way as for algebras, it is possible to define a tensor 
product between modules (cf. \cite{Mit2}):
\begin{defi}
Let $M$ be a left $\C$-module and let $N$ be a right $\C$-module. 
The \textbf {tensor product over $\C$ between $M$ and $N$}, $M 
\otimes_{\C} N$, is defined as the $k$-module given by 
\[     M \otimes_{\C} N = (\bigoplus_{x \in \C_{0}} M_{x} \otimes 
{}_{x}N) / 
       \cl{\{ m.f \otimes n - m \otimes f.n : m \in M_{x}, n \in 
{}_{y}N, f \in {}_{x}\C_{y} \}}.     \]
If $M$ and $N$ are $\C$-bimodules, it is also possible to define the 
$\C$-\textbf {bimodule tensor product over $\C$}:  
\[     {}_{y}(M \otimes_{\C} N)_{x} = (\bigoplus_{z \in \C_{0}} 
{}_{y}M_{z} \otimes_{k} {}_{z}N_{x}) / 
       \cl{\{ m.f \otimes n - m \otimes f.n \}},     \]
where $m \in {}_{y}M_{y'}$, $n \in {}_{x'}N_{x}$, $f \in 
{}_{y'}\C_{x'}$. 
\end{defi}

Next we recall the definition of Hochschild-Mitchell homology and 
cohomology. 
Standard (co)homological methods are available in ${}_{\C}\mathrm{Mod}$ 
(cf. \cite{HS1}). 

\begin{defi}
\label{obstor}
Let $( x_{n+1},\dots,x_{1} )$ be a $(n+1)$-sequence of objects of $\C$. 
The \textbf{$k$-nerve associated to the $(n+1)$-sequence} is the 
$k$-vector space 
\[     {}_{x_{n+1}}\C_{x_{n}} \otimes \dots \otimes 
{}_{x_{2}}\C_{x_{1}}.     \]
The \textbf{$k$-nerve of $\C$ in degree $n$} ($n \in \NN_{0}$)  is 
\[     \overline {N}_{n} = \bigoplus_{(n+1)\hbox{\footnotesize-tuples}} 
{}_{x_{n+1}}\C_{x_{n}} \otimes \dots 
        \otimes {}_{x_{2}}\C_{x_{1}}.     \]      

There is $\C^{e}$-bimodule associated to $\overline{N}_{n}$ defined by  
\[     {}_{y}(N_{n})_{x} = \bigoplus_{(n+1)\hbox{\footnotesize-tuples}} 
{}_{y}\C_{x_{n+1}} \otimes {}_{x_{n+1}}\C_{x_{n}} 
       \otimes \dots \otimes {}_{x_{2}}\C_{x_{1}} \otimes 
{}_{x_{1}}\C_{x}.     \]    

Then the associated \textbf{Hochschild-Mitchell complex} is
\[     \dots \stackrel{d_{n+1}}{\rightarrow} N_{n} 
\stackrel{d_{n}}{\rightarrow} \dots
         \stackrel{d_{2}}{\rightarrow} N_{1} 
\stackrel{d_{1}}{\rightarrow}
         N_{0} \stackrel{d_{0}}{\rightarrow} \C \rightarrow 0,     \] 
where $d_{n}$ is given by the usual formula, i.e. 
\[     d_{n} (f_{0} \otimes \dots \otimes f_{n+1}) = \sum\limits_{k = 
0}^{n} (-1)^{k} f_{0} \otimes \dots \otimes f_{k}.f_{k+1} \otimes \dots 
\otimes f_{n+1}.     \]

This complex is a projective resolution of the $\C$-bimodule $\C$. 
The proof that it is a resolution is similar to the standard proof for 
algebras. 
\end{defi}

\begin{defi}
Given a $\C$-bimodule $M$ the \textbf{Hochschild-Mitchell cohomology of 
$\C$ with coefficients in $M$} is the cohomology of the following 
cochain 
complex 
\[     0 \rightarrow \prod_{x \in \C_{0}} {}_{x}M_{x} 
\stackrel{d^{0}}{\rightarrow} 
        \mathrm{Hom}(N_{1},M) \stackrel{d^{1}}{\rightarrow} \dots 
\stackrel{d^{n-1}}{\rightarrow} 
       \mathrm{Hom}(N_{n},M) \stackrel{d^{n}}{\rightarrow} \dots,     
\]
where $d$ is given by the usual formula, and   
  \begin{align*}
    C^{n}(\C , M) &= \mathrm{Hom}(N_{n} , M) = \mathrm{Nat}(N_{n} , M)      
    \\
    &= \prod_{(n+1)\hbox{\footnotesize-tuples}} \mathrm{Hom}_{k} 
({}_{x_{n+1}}\C_{x_{n}} \otimes \dots \otimes 
       {}_{x_{2}}\C_{x_{1}},{}_{x_{n+1}}M_{x_{1}}).    
  \end{align*} 
We denote it $H^{\bullet}(\C , M)$. 

Analogously the \textbf{Hochschild-Mitchell homology of $\C$ with 
coefficients 
in $M$} is the homology of the chain complex 
\[     \dots \stackrel{d_{n+1}}{\rightarrow} M \otimes N_{n} 
\stackrel{d_{n}}{\rightarrow} \dots
         \stackrel{d_{2}}{\rightarrow} M \otimes N_{1} 
\stackrel{d_{1}}{\rightarrow}
       \bigoplus_{x \in \C_{0}} {}_{x}M_{x} \rightarrow 0,     \] 
where $d$ is given by the usual formula and
  \[
    C_{n}(\C , M) = M \otimes_{\C^{e}} N_{n}  
    = \bigoplus_{(n+1)\hbox{\footnotesize-tuples}} 
    {}_{x_{1}}M_{x_{n+1}} \otimes {}_{x_{n+1}}\C_{x_{n}} \otimes \dots 
\otimes {}_{x_{2}}\C_{x_{1}}.     
  \] 
We denote it $H_{\bullet}(\C , M)$. 
\end{defi}

The following is a generalization of Watt's Theorem for modules over 
$k$-algebras:
\begin{teo}
\label{thm:1}
Let $\C$ and $\D$ be $k$-linear categories and let $F : 
{}_{\C}\mathrm{Mod} \rightarrow {}_{\D}\mathrm{Mod}$ be a functor. 
The following statements are equivalent 
\begin{itemize}
\item[(a)]
$F$ preserves arbitrary direct sums and is right exact.

\item[(b)]
There exists a $\D$-$\C$-bimodule $T$ such that $F(\place) = T 
\otimes_{\C} (\place)$.

\item[(c)]
$F$ has a right adjoint.
\end{itemize}
\end{teo}
\noindent\textbf{Proof.}
We trivially have that (b) implies (c), and (c) implies (a). 
Let us prove that (a) implies (b). 
For each $x \in \C_{0}$ define the left $\D$-module 
\[     {}_{\place}T_{x} = F( {}_{\place}\C_{x} ).     \]
The collection $\{ {}_{y}T_{x} \}_{y \in \D_{0},x\in \C_{0}}$ is a 
$\D$-$\C$-bimodule as we shall now prove: 
it is trivially a left $\D$-module by definition. 
Given ${}_{x}f_{x'} \in {}_{x}\C_{x'}$, it induces a morphism of left 
$\C$-modules 
\[     .{}_{x}f_{x'} : {}_{\place}\C_{x} \rightarrow 
{}_{\place}\C_{x'},     \]
given by right multiplication by ${}_{x}f_{x'}$, so we get a morphism 
of left $\D$-modules 
\[     F(.{}_{x}f_{x'}) : F({}_{\place}\C_{x}) \rightarrow 
F({}_{\place}\C_{x'}).     \]
This natural transformation gives the structure of right $\C$-module. 

Moreover, both actions are compatible since the map $F(.{}_{x}f_{x'})$ 
is a morphism of left $\D$-modules. 

We have that 
\[     {}_{y}T \otimes_{\C} \C_{x} 
= (\oplus_{z \in \C_{0}} {}_{y}T_{z} \otimes {}_{z}\C_{x})/\cl{\{ t_{z} 
{}_{z}f_{z'} \otimes {}_{z'}g - t_{z} \otimes {}_{z}f_{z'} {}_{z'}g\}} 
\simeq {}_{y}T_{x},     \]
using the $k$-linear isomorphism $\overline{t \otimes f} \mapsto tf$ 
(with inverse $t \mapsto \overline{t \otimes 1}$). 
This gives naturally a left $\D$-module isomorphism 
\[     {}_{\place}T \otimes_{\C} \C_{x} = {}_{\place}T_{x}.     \]

From now on, given $M \in {}_{\C}\mathrm{Mod}$, we are going to write 
$F(M)$ instead of $F({}_{\place}M)$. 
We shall now prove that $F(\place) = T \otimes_{\C} (\place)$. 
Since $F$ and $\otimes_{\C}$ commute with direct sums, there are 
isomorphisms of $\D$-modules 
\[     F(\oplus_{i\in I} \C_{x_{i}}) = \oplus_{i\in I} F(\C_{x_{i}}) = 
\oplus_{i \in I} T_{x_{i}} = \oplus_{i \in I} T \otimes_{\C} \C_{x_{i}} 
= T \otimes_{\C} (\oplus_{i\in I} \C_{x_{i}}).     \] 

Given any left $\C$-module $M$ there is an exact sequence 
\[     \oplus_{j \in J} \C_{y_{j}} \rightarrow \oplus_{i \in I} 
\C_{x_{i}} \rightarrow M \rightarrow 0,     \]
hence, by right exactness, we get that
\[     F(\oplus_{j \in J} \C_{y_{j}}) \rightarrow F(\oplus_{i\in I} 
\C_{x_{i}}) \rightarrow F(M) \rightarrow 0     \]
is exact. 
Taking into account the previous isomorphism, the following diagram has 
exact rows and commuting squares   
\[
\xymatrix
{
  F(\bigoplus_{j \in J} \C_{y_{j}}) 
  \ar[r]
  \ar@{=}[d]
  &  
  F(\bigoplus_{i \in I} \C_{x_{i}}) 
  \ar[r]
   \ar@{=}[d]
  &
  F(M) 
  \ar[r]
  &
  0
  \\
  \bigoplus_{j \in J} T \otimes_{\C} \C_{y_{j}} 
  \ar[r]
  &  
  \bigoplus_{i \in I} T \otimes_{\C} \C_{x_{i}} 
  \ar[r]
  &
  T \otimes_{\C} M 
  \ar[r]
  &
  0.
}
\]

By diagrammatic considerations we get a map $F(M) \rightarrow T 
\otimes_{\C} M$ making the whole diagram commutative. 
Then the Five lemma assures that this map is an isomorphism. 

The naturality of the map is also clear: 
if $f : M \rightarrow N$ is a $\C$-module morphism, there is a 
commutative diagram 
\[
\xymatrix
{
  \bigoplus_{j \in J} \C_{y_{j}} 
  \ar[r]
  \ar[d]
  &  
  \bigoplus_{i \in I} \C_{x_{i}} 
  \ar[r]
   \ar[d]
  &
  M 
  \ar[r]
  \ar[d]^{f}
  &
  0
  \\
  \bigoplus_{j' \in J'} \C_{y'_{j'}}
  \ar[r]
  &  
  \bigoplus_{i' \in I'} \C_{x'_{i'}} 
  \ar[r]
  &
  N 
  \ar[r]
  &
  0
}
\]

and hence a diagram 
\[
\xymatrix@C-20pt
{
  F(\bigoplus\limits_{j \in J} \C_{y_{j}}) 
  \ar[rr]
  \ar[rd]
  \ar@{=}[dd]
  &
  &  
  F(\bigoplus\limits_{i \in I} \C_{x_{i}}) 
  \ar[rr]
  \ar[rd]
  \ar@{=}[dd]
  &
  &
  F(M) 
  \ar[rr]
  \ar[rd]
  \ar[dd]
  &
  &
  0
  &
  \\
  &
  F(\bigoplus\limits_{j' \in J'} \C_{y'_{j'}}) 
  \ar[rr]
  \ar@{=}[dd]
  &
  &  
  F(\bigoplus\limits_{i' \in I'} \C_{x'_{i'}}) 
  \ar[rr]
  \ar@{=}[dd]
  &
  &
  F(N) 
  \ar[rr]
  \ar[dd]
  &
  &
  0
  \\
  \bigoplus\limits_{j \in J} T \otimes_{\C} \C_{y_{j}} 
  \ar[rr]
  \ar[rd]
  &
  &  
  \bigoplus\limits_{i\in I} T \otimes_{\C} \C_{x_{i}} 
  \ar[rr]
  \ar[rd]
  &
  &
  T \otimes_{\C} M 
  \ar[rr]
  \ar[rd]
  &
  &
  0
  &
  \\
  &
  \bigoplus\limits_{j' \in J'} T \otimes_{\C} \C_{y'_{j'}} 
  \ar[rr]
  &
  &  
  \bigoplus\limits_{i' \in I'} T \otimes_{\C} \C_{x'_{i'}} 
  \ar[rr]
  &
  &
  T \otimes_{\C} N
  \ar[rr]
  &
  &
  0
}
\]
Since the two left vertical faces (normal to the page) commute (by 
definition of $T$), 
we obtain that the right vertical face also commutes, and this fact 
proves the naturality. 
\qed

As an application of the characterization of such functors we obtain a 
description of Morita equivalences 
of $k$-linear categories. 
We also give an example relating this description to the one given in 
\cite{CS1}. 

\begin{teo}
\label{thm:2}
Let $\C$ and $\D$ be two $k$-linear categories. 
They are (left) Morita equivalent if and only if there are a 
$\C$-$\D$-bimodule $P$ and 
a $\D$-$\C$-bimodule $Q$ such that $P \otimes_{\D} Q \simeq \C$ and $Q 
\otimes_{\C} P \simeq \D$ as bimodules. 
Furthermore, these bimodules satisfy that $\{P_{y}\}_{y \in \D_{0}}$ 
and $\{ Q_{x} \}_{x \in \C_{0}}$ are sets of 
projective and finitely generated generators of ${}_{\C}\mathrm{Mod}$ 
and ${}_{\D}\mathrm{Mod}$ respectively. 
\end{teo}
\noindent\textbf{Proof.}
Given two bimodules $P$ and $Q$, we define the functors 
\begin{align*}
   F : {}_{\C}\mathrm{Mod} &\rightarrow {}_{\D}\mathrm{Mod}
   \\
   F(\place) &= Q \otimes_{\C} (\place),
\end{align*}
and 
\begin{align*}
   G : {}_{\D}\mathrm{Mod} &\rightarrow {}_{\C}\mathrm{Mod}
   \\
   G(\place) &= P \otimes_{\D} (\place).  
\end{align*}
Since $P \otimes Q$ and $Q \otimes P$ are isomorphic as bimodules to 
$\C$ and $\D$ respectively, 
then $F \circ G \simeq \id_{\D}$ and $G \circ F \simeq \id_{\C}$. 

Conversely, let $F : {}_{\C}\mathrm{Mod} \rightarrow 
{}_{\D}\mathrm{Mod}$ be a functor giving the equivalence 
with quasi-inverse functor $G$. 
Since an equivalence preserves direct sums and is exact, Theorem 
\ref{thm:1} guarantees the existence of a $\C$-$\D$-bimodule $P$ and a 
$\D$-$\C$-bimodule $Q$ satisfying 
\begin{align*}
   F(\place) &= Q \otimes_{\C} (\place)
   \\
   G(\place) &= P \otimes_{\D} (\place).  
\end{align*}

The isomorphism $F \circ G \simeq \id_{\D}$ implies that 
$\D \simeq F \circ G (\D) = Q \otimes_{\C} (P \otimes_{\D} \D) \simeq Q 
\otimes_{\C} P$. 
The other isomorphism is analogous. 

Since, given $x \in \C_{0}$, $Q_{x}$ is isomorphic to $F(\C_{x})$, each 
$Q_{x}$ is finitely generated and projective, 
and the same applies to $P_{y}$ ($y \in \D_{0}$). 
Also, taking into account that $\{ \C_{x} \}_{x \in \C_{0}}$ is a set 
of generators of ${}_{\C}\mathrm{Mod}$ and $F$ is an 
equivalence, we get that $\{ Q_{x} \}_{x \in \C_{0}} = \{ F(\C_{x}) 
\}_{x \in \C_{0}}$ is a set of generators of 
${}_{\D}\mathrm{Mod}$. 
The same arguments apply to $\{ P_{y} \}_{y \in \C_{0}}$. 
\qed

\begin{rem}
We infer from the theorem above that if $\C$ and $\D$ are left Morita 
equivalent, then they are right Morita 
equivalent.  
This is done just by taking the functors 
\begin{align*}
   F : \mathrm{Mod}_{\C} &\rightarrow \mathrm{Mod}_{\D}
   \\
   F(\place) &= (\place) \otimes_{\C} P,
\end{align*}
and 
\begin{align*}
   G : \mathrm{Mod}_{\D} &\rightarrow \mathrm{Mod}_{\C}
   \\
   G(\place) &= (\place) \otimes_{\D} Q.  
\end{align*}
\end{rem}

The following results will complete the description. 
\begin{prop}
Let $\C$, $\D$, $\E$ and $\F$ be $k$-linear categories and 
${}_{\D}P_{\E}$, ${}_{\D}M_{\C}$, ${}_{\C}N_{\F}$ be a set 
of bimodules. 
Then the following is a natural morphism of $\E$-$\F$-bimodules   
\[    \eta : \mathrm{Hom}_{\D} (P , M) \otimes_{\C} N \rightarrow 
\mathrm{Hom}_{\D} (P , M \otimes_{\C} N),     \]
defined by 
\[     \eta(t \otimes n)(p) = t(p) \otimes n.     \]

Furthermore, if $P_{x}$ is finitely generated and projective as left 
$\D$-module for each $x \in \E_{0}$, 
then $\eta$ is an isomorphism. 
\end{prop}
\noindent\textbf{Proof.}
The morphism is clearly well-defined and natural. 
To prove the second statement, let us first suppose that $\E = \D$ and 
$P = {}_{\place}\D_{\place}$. 
Since for each $x \in \D_{0}$ we have an isomorphism of right 
$\D$-modules 
\[     \mu_{M} : \mathrm{Hom}_{\D} (\D_{x} , M) 
\overset{\simeq}{\rightarrow} {}_{x}M     \]
defined via the Yoneda's isomorphism 
\[     \mu_{M}(t) = {}_{x}t_{x} (\id_{x}),     \]
we get 
\[     \mu_{M} \otimes \id: \mathrm{Hom}_{\D} (\D_{x} , M) \otimes_{\C} 
N 
       \overset{\simeq}{\rightarrow} {}_{x}M \otimes_{\C} N,     \]
and 
\[     \mu_{M \otimes_{\C} N} : \mathrm{Hom}_{\D} (\D_{x} , M 
\otimes_{\C} N) 
       \overset{\simeq}{\rightarrow} {}_{x}M \otimes_{\C} N.     \]
We see that $\eta = \mu_{M \otimes_{\C} N}^{-1} \circ (\mu_{M} \otimes 
\id)$. 

Now if, $P_{x}$ is finitely generated and projective, there exists $P'$ 
such that 
\[     P' \oplus P_{x} = \bigoplus_{i = 1}^{n} \C_{x_{i}}.     \]
Using lemma (20.9) from \cite{AF1}, we are able to prove that $\eta$ is 
an isomorphism.  
\qed

The proof of the following proposition is analogous:
\begin{prop}
Let $\C$, $\D$, $\E$ and $\F$ be $k$-linear categories and 
${}_{\E}P_{\D}$, ${}_{\C}M_{\D}$, ${}_{\C}N_{\F}$ be a set 
of bimodules. 
Then the following is a natural morphism of $\E$-$\F$-bimodules  
\begin{align*}
    \nu : P \otimes_{\D} \mathrm{Hom}_{\C} (M , N) 
&\overset{\simeq}{\rightarrow} \mathrm{Hom}_{\C} (\mathrm{Hom}_{\D} (P , M) , N),     
    \\
    \nu(p \otimes t)(u) &= t(u(p)).     
\end{align*}

Furthermore, if ${}_{x}P$ is finitely generated and projective as right 
$\D$-module for each $x \in \E_{0}$, 
then $\nu$ is an isomorphism. 
\end{prop}

From the previous propositions we obtain: 
\begin{coro}
Two $k$-linear categories $\C$ and $\D$ are Morita equivalent if and 
only if there exists a $\C$-$\D$-bimodule $P$ 
such that $\{ P_{y} \}_{y \in \D_{0}}$ is a set of finitely generated 
projective generators of ${}_{\C}\mathrm{Mod}$, 
$\{ {}_{x}P \}_{x \in \C_{0}}$ is a set of finitely generated 
projective generators of $\mathrm{Mod}_{\D}$ 
and $Hom_{\C} (P,P) \simeq \D$ (as $\D$-bimodules).
\end{coro}
\noindent\textbf{Proof.}
If $\C$ and $\D$ are Morita equivalent then we use the bimodule $P$ 
defined in Theorem \ref{thm:2} which satisfies 
all the conditions except perhaps that $\mathrm{Hom}_{\C}(P,P) = \D$. 
But $Q \otimes_{\C} \place = \mathrm{Hom}_{\C} (P, \place)$, so we get that 
$\D = Q \otimes_{\C} P = \mathrm{Hom}_{\C} (P, P)$. 

Conversely, suppose that there exists a $\C$-$\D$-bimodule $P$ such 
that $\{ P_{y} \}_{y \in \D_{0}}$ is 
a set of finitely generated projective generators of 
${}_{\C}\mathrm{Mod}$ and $Hom_{\C} (P,P) \simeq \D$ as $\D$-bimodules. 
Then we set 
\begin{align*}
   F : {}_{\C}\mathrm{Mod} &\rightarrow {}_{\D}\mathrm{Mod}
   \\
   F(\place) &= \mathrm{Hom}_{\C} (P,\place),
\end{align*}
and 
\begin{align*}
   G : {}_{\D}\mathrm{Mod} &\rightarrow {}_{\C}\mathrm{Mod}
   \\
   G(\place) &= P \otimes_{\D} (\place).  
\end{align*}

From the previous proposition we have that, for any left $\D$-module 
$M$, 
\[     \mathrm{Hom}_{\C} (P , P \otimes_{\D} M) \simeq 
\mathrm{Hom}_{\C} (P , P) \otimes_{\D} M 
       \simeq \D \otimes_{\D} M \simeq M,     \]
and also, for any left $\C$-module $N$, 
\[     P \otimes_{\D} \mathrm{Hom}_{\C} (P , N) \simeq 
\mathrm{Hom}_{\D} (\mathrm{Hom}_{\C}(P , P) , N)  
       \simeq \mathrm{Hom}_{\D} (\D , N) \simeq N,     \]
where all isomorphisms are natural. 
Hence $F$ and $G$ are quasi-inverse functors, giving the Morita 
equivalence. 
\qed

\begin{ejem}
Suppose that $E$ is a partition of the set of objects $\C_{0}$ of a 
$k$-linear category $\C$, given by 
$\C_{0} = \sqcup_{e \in E} E_{e}$, with $\#(E_{e}) < \aleph_{0}$, 
$\forall e \in E$. 
It is proved in \cite{CS1} that $\C$ is Morita equivalent to the 
contracted category along the partition $E$, $\C/E$. 
In fact the functors giving the equivalence are the following 
\begin{align*}
   F : {}_{\C}\mathrm{Mod} &\rightarrow {}_{\C/E}\mathrm{Mod} 
   \\
   {}_{e}F(M) &= \bigoplus_{x \in E_{e}} {}_{x}M,
\end{align*}
and
\begin{align*}
G : {}_{\C/E}\mathrm{Mod} &\rightarrow {}_{\C}\mathrm{Mod} 
\\
{}_{x}G(N) &= f_{x}.({}_{e}N),
\end{align*}
where $e$ is the unique element of $E$ such that $x \in E_{e}$ and 
$f_{x}$ is the idempotent 
$\mid E\mid \times \mid E\mid$-matrix. 

The bimodules giving the equivalence are: 
\begin{align*}
{}_{e}({}_{\C/E}P_{\C})_{x} &= \bigoplus_{y \in E_{e}} {}_{y}\C_{x},
\\
{}_{x}({}_{\C}Q_{\C/E})_{e} &= \bigoplus_{y \in E_{e}} {}_{x}\C_{y}.
\end{align*}
It is easy to check that $P \otimes_{\C} Q \simeq \C/E$, $Q 
\otimes_{\C/E} P \simeq \C$ as bimodules and also 
$F(\place) = P \otimes_{\C} (\place)$ and $G(\place) = Q \otimes_{\C/E} 
(\place)$.
\end{ejem}

From now on we shall consider the derived category of 
${}_{\C}\mathrm{Mod}$. 
This is a special case of the theory developed by Keller for DG 
categories. 
We will recall some definitions, but we refer the reader to 
\cite{Kel1}, for further references.  

As usual, we consider the category of $\C$-modules embedded into the 
category of cochains of complexes of $\C$-modules, 
denoted by $\mathrm{Ch} ({}_{\C}\mathrm{Mod})$ or $\mathrm{Ch} (\C)$, 
and we denote the shift of a complex 
$M^{\bullet}$ by $M^{\bullet}[1]$ or $SM^{\bullet}$, the homotopy 
category by $\H({}_{\C}\mathrm{Mod})$ 
or just by $\H(\C)$, and the derived category by 
$D({}_{\C}\mathrm{Mod})$ or $D(\C)$. 

We say that a complex of $\C$-modules $M$ is \textbf{relatively 
projective} if it is a direct summand of a direct sum of complexes of the form 
$\C_{x}[n]$, for $n \in \ZZ$, $x \in \C_{0}$. 
We also recall that a complex of $\C$-modules $M$ is 
\textbf{homotopically projective} if it is homotopically equivalent 
to a complex $P$ provided with an increasing filtration (indexed by 
$\NN_{0}$)
\[     P_{-1} = 0 \subset P_{0} \subset \dots \subset P_{n} \subset 
\dots \subset P,     \]
satisfying the following properties: 
\begin{enumerate}
\item $P = \bigcup_{n \in \NN_{0}} P_{n}$. 

\item The inclusion $P_{n} \subset P_{n+1}$ ($n \in \NN_{0}$) splits in 
the category of graded modules over $\C$. 

\item The quotient $P_{n}/P_{n-1}$ ($n \in \NN_{0}$) is isomorphic in 
$\mathrm{Ch}(\C)$ to a relatively projective 
module.  
\end{enumerate}
As it is proved in \cite{Kel1}, the following is a split exact sequence 
in the category of graded modules over $\C$ 
\begin{equation}
\label{eq:seq}
     \bigoplus_{n \in \NN_{0}} P_{n} \rightarrow \bigoplus_{n \in 
\NN_{0}} P_{n} \rightarrow P,     
\end{equation}
and this split exact sequence gives a triangle in $\H(\C)$. 

We denote $\H_{p}(\C)$ the full triangulated subcategory of $\H(\C)$ 
formed by homotopically projective complexes of modules. 
We denote $\H_{p}^{b}(\C)$ the smallest strictly (i.e., closed under 
isomorphisms) full triangulated subcategory of $\H_{p}(\C)$ containing 
$\{ \C_{x} \}_{x \in \C_{0}}$. 

We recall the following theorem from \cite{Kel1}, pp. 69--70, Thm. 3.1: 
\begin{teo}
For any complex of $\C$-modules $M$ we have the following triangle in 
$\mathcal{H}{\C}$ 
\[     p(M) \rightarrow M \rightarrow M \rightarrow a(M) \rightarrow S 
p(M),     \]
where $a(M)$ is acyclic and $p(M)$ is homotopically projective. 

Furthermore, this construction gives rise to triangle functors $p$ and 
$a$ on $\mathcal{H}(\C)$ commuting 
with direct sums, $p$ is the right adjoint of the inclusion functor 
from the full triangulated subcategory 
of homotopically projective complexes, and $a$ is the left adjoint of 
the inclusion of the full triangulated 
subcategory of acyclic complexes. 
\end{teo}
Following Keller, we call $p(M)$ \textbf{the projective resolution} of 
the complex $M$. 

Taking into account that any $k$-linear category is a DG category 
concentrated in degree $0$ with null differential, 
we may apply the following theorem (cf. \cite{Kel1}, corollary 9.2), 
adapted to the $k$-linear case, 
\begin{teo}
\label{thm:Kel9.2}
Let $\C$ and $\D$ be two $k$-linear categories such that $\D$ is 
$k$-flat 
(i.e., ${}_{y}\D_{x}$ is $k$-flat, for every $x, y \in \D_{0}$). 
The following are equivalent: 
\begin{itemize}
\item[(i)]
There is a $\C$-$\D$-bimodule $P$ such that $P \otimes^{L}_{\C} \place 
: D(\C) \rightarrow D(\D)$ is an equivalence. 

\item[(ii)] 
There is an $S$-equivalence $D(\C) \rightarrow D(\D)$. 

\item[(iii)]
$\C$ is equivalent to a full subcategory $\E$ of $D(\D)$ whose objects 
form a set of small generators and 
satisfy the following 
\[     \mathrm{Hom}_{D(\D)} (M, N[n]) = 0,     \]
for all $n \neq 0$, $M, N \in \E$. 
\end{itemize}
\end{teo}
If any of the three equivalent conditions of the theorem is satisfied 
we say that $\C$ and $\D$ are derived equivalent. 

We recall that a $k$-linear category is \textbf{projective} if 
${}_{y}\C_{x}$ is a projective $k$-module for every $x, y \in \C_{0}$. 
We obtain the following as a corollary of the previous theorem.  
In particular, the hypothesis of projectivity holds when $k$ is a 
field. 
\begin{teo}
\label{thm:hoch}
Let $\C$ and $\D$ be two small $k$-linear projective categories which 
are derived equivalent. 
Then the Hochschild-Mitchell homology and cohomology groups of $\C$ and 
$\D$ are respectively isomorphic. 
\end{teo}
\noindent\textbf{Proof.}
Since $\C$ and $\D$ are derived equivalent there exists a 
$\D$-$\C$-bimodule $P$ and a $\C$-$\D$-bimodule $Q$, 
such that 
\[     P \otimes^{L}_{\C} \place \otimes^{L}_{\C} Q : D(\C^{e}) 
\rightarrow D(\D^{e}),     \]
is an equivalence, with quasi-inverse 
\[     Q \otimes^{L}_{\C} \place \otimes^{L}_{\C} P : D(\D^{e}) 
\rightarrow D(\C^{e}).     \]

As a consequence, $P \otimes^{L}_{\C} Q \simeq \D$ in $D(\D^{e})$, and 
$Q \otimes^{L}_{\D} P \simeq \C$ in $D(\C^{e})$. 

Hence, we have the following chain of isomorphisms in $D(k)$ 
\[     \C \otimes^{L}_{\C^{e}} \C \overset{\simeq}{\rightarrow}  
       (Q \otimes_{\D}^{L} P) \otimes^{L}_{\C^{e}} (Q \otimes_{\D}^{L} 
P) \overset{\simeq}{\rightarrow} 
       (P \otimes_{\C}^{L} Q) \otimes^{L}_{\D^{e}} (P \otimes_{\C}^{L} 
Q) \overset{\simeq}{\rightarrow} 
       \D \otimes^{L}_{\D^{e}} \D,     \]
where the second isomorphism is induced by the isomorphism 
\begin{align*}
  (p(Q) \otimes_{\D} p(P)) \otimes_{\C^{e}} (p(Q) \otimes_{\D} p(P)) 
&\rightarrow 
  (p(P) \otimes_{\C} p(Q)) \otimes_{\D^{e}} (p(P) \otimes_{\C} p(Q))     
  \\
 (a \otimes b) \otimes (a' \otimes b') &\mapsto (b \otimes a') \otimes 
(b' \otimes a),   
\end{align*}  
and the fact that $p(P) \otimes_{\C} p(Q)$ is a projective resolution 
of $P \otimes_{\C} Q$ in ${}_{\D^{e}}\mathrm{Mod}$ 
and $p(Q) \otimes_{\D} p(P)$ is a projective resolution of $Q 
\otimes_{\C} P$ in ${}_{\C^{e}}\mathrm{Mod}$. 
To prove this last statement we proceed as follows. 
Since $\C$ and $\D$ are projective $k$-categories, given a 
homotopically projective $\C$-$\D$-bimodule $M$ 
(which we may suppose of the form $(\C_{x} \otimes_{k} {}_{y}\D)$ for 
$x \in \C_{0}$, $y \in \D_{0}$) 
the functor $M \otimes_{\D} \place$ sends relatively projective 
$\D$-$\C$-bimodules of type 
$(\D_{y'} \otimes_{k} {}_{x'}\C)$ (for $x' \in \C_{0}$, $y' \in 
\D_{0}$) into 
$\C_{x} \otimes_{k} {}_{y}\D_{y'} \otimes_{k} {}_{x'}\C$, which are 
relatively projective $\C$-bimodules. 
Hence we get that $M \otimes_{\D} \place$ sends homotopically 
projective $\D$-$\C$-bimodules into homotopically projective $\C$-bimodules. 

This implies immediately the theorem for homology, since we have 
\[     H_{n} (\C \otimes_{\C^{e}}^{L} \C) \simeq 
\mathrm{Tor}_{n}^{\C^{e}} (\C,\C) = HH_{n} (\C).     \]

For cohomology, we make use of the following isomorphism 
\[     \mathrm{Hom}_{D(\C^{e})} (\C, \C[n]) \simeq 
\mathrm{Ext}^{n}_{\C^{e}} (\C, \C) = HH^{n} (\C).     \]
which is proved in the second lemma of \cite{Kra1} section 1.5. 
This concludes the proof of the theorem. 
\qed

\section{Derived equivalences between one-point extensions}

Let us first state some facts concerning convex categories. 
\begin{defi}
Let $\C$ be a linear $k$-category and $\D$ a subcategory. 
We say that $\D$ is a \textbf{convex subcategory of $\C$} if given 
$x_{0}, x_{n} \in \D_{0}$, 
$x_{1} , \dots, x_{n-1} \in \C_{0}$ such that $\exists i$, $1 \leq i 
\leq n-1$ with $x_{i} \notin \D_{0}$, 
and morphisms $f_{i} \in {}_{x_{i+1}}\C_{x_{i}}$, for $0 \leq i \leq 
n-1$ 
then $f_{n-1} \circ \dots \circ f_{0} = 0$. 
\end{defi}

\begin{obs}
The following facts about convex categories are easy to prove:
\begin{itemize}
\item If $\C'$ is a convex subcategory of $\C$, then $\C'^{op}$ is a 
convex subcategory of $\C^{op}$. 

\item If $\C'$ and $\D'$ are convex subcategories of $\C$ and $\D$ respectively, 
then $\C' \boxtimes \D'$ is a convex subcategory of $\C \boxtimes \D$. 
\end{itemize}
\end{obs}

If $\D$ is a convex subcategory of $\C$ then there is a functor 
\[     i : \mathrm{Mod}_{\D} \rightarrow \mathrm{Mod}_{\C}     \]
given by the $i(N)_{x} = N_{x}$, for $x \in \D_{0}$ and $i(N)_{y} = 0$, 
for $y \in \C_{0} \setminus \D_{0}$. 
The action of $\C$ is induced by the action of $\D$ on $N$. 
It is clear that $i(N)$ is a right $\C$-module and it is well defined 
since $\D \subset \C$ is convex. 

There is a functor induced by the inclusion $r : \mathrm{Mod}_{\C} \rightarrow \mathrm{Mod}_{\D}$, 
given by $r(M)= M \circ inc_{\D \subset \C}$. 

They are adjoint functors, namely, we have the isomorphism
\begin{align*} 
    \theta : \mathrm{Hom}_{\D} (r(M) , N) &\rightarrow 
\mathrm{Hom}_{\C} (M , i(N))  
\\
    \theta (\{t_{y}\}_{y \in \D_{0}})_{x} &= \begin{cases} t_{x} \hskip 0.5cm \text{if $x \in D_{0}$}
                                                             \\
                                                             0 \hskip 0.5cm \text{otherwise.}
                                              \end{cases}
\end{align*}
It is easy to check that this map is well defined and it is natural, 
and it is an isomorphism with inverse 
is given by
\begin{align*} 
    \zeta : \mathrm{Hom}_{\C} (M , i(N)) &\rightarrow \mathrm{Hom}_{\D} 
(r(M) , N)  
\\
    \zeta (\{t_{x}\}_{x \in \C_{0}})_{y} &= t_{y}, \hskip 0.5cm 
\text{for $y \in \D_{0}$.}                    
\end{align*}

The adjunction says immediately that $r$ preserves epimorphisms and $i$ 
preserves monomorphisms, 
but we may also easily see that $r$ preserves monomorphisms and $i$ 
preserves epimorphisms. 
Hence both functors are exact, $r$ preserves projective objects and $i$ 
preserves injective objects. 

\begin{lema}
\label{lema:convex}
Let $\D$ be a full convex subcategory of $\C$ and let $M, N$ be two 
$\D$-modules. 
Then there is an isomorphism 
\[     \mathrm{Ext}^{\bullet}_{\D} (M , N) \simeq 
\mathrm{Ext}_{\C}^{\bullet} (i(M) , i(N)).     \]
\end{lema}
\noindent\textbf{Proof.}
Choosing a projective $\C$-resolution $P_{\bullet}$ of $i(M)$, since 
$r$ is exact and 
preserves projectives, then $r(P_{\bullet})$ is a projective $\D$-resolution 
of 
$r(i(M)) = M$. 
By the previous adjunction there is a morphism of complexes 
\[     \mathrm{Hom}_{\C} (P_{\bullet} , i(N)) \simeq \mathrm{Hom}_{\D} 
(r(P_{\bullet}) , N),     \]
implying that 
\[     \mathrm{Ext}^{\bullet}_{\D} (M , N) \simeq 
\mathrm{Ext}_{\C}^{\bullet} (i(M) , i(N)).     \] 
\qed  

Next let us define, given a $k$-linear category $\C$ and a right 
$\C$-module $M$, \textbf{the one point extension of $\C$ by $M$} as the 
following small category, which we will denote $\C[M]$. 
The set of objects is $(\C[M])_{0} = \C_{0} \sqcup \{ M \}$. 
The set of morphisms is given by 
\[ {}_{y}\C[M]_{x} = {}_{y}\C_{x},  
\hskip 0.5cm
{}_{M}\C[M]_{x} = M_{x}, 
\hskip 0.5cm
{}_{y}\C[M]_{M} = 0, 
\hskip 0.5cm
{}_{M}\C[M]_{M} = k,
\hskip 0.5cm
\text{for $x, y \in \C_{0}$.} 
\]
The composition is given by composition in $\C$, the action of $\C$ on 
$M$ 
and the structure of $k$-module on each $M_{x}$. 
It may be easily verified that $\C[M]$ satisfies the axioms of a 
$k$-linear category and 
that $\C$ is a full convex subcategory of $\C[M]$. 

\begin{obs}
There is a dual definition for a left $\C$-module $M$, the only changes 
are ${}_{x}\C[M]_{M} = {}_{x}M$ and 
${}_{M}\C[M]_{y} = 0$. 

If $\C$ is finite, then $a(\C[M]) \simeq a(\C)[M]$, where the last one 
denotes the one point extension 
of the algebra $a(\C)$ by the induced module $\oplus_{x \in \C_{0}} 
M_{x}$. 
\end{obs}

In this context, we define the right $\C[M]$-module $\bar{M}$, given by 
$\bar{M}_{x} = M_{x}$ ($x \in \C_{0}$) and $\bar{M}_{M} = k$. 
The action is the following 
\[
\xymatrix@R-20pt
{
\bar{\rho}_{x,y} : \bar{M}_{x} \otimes {}_{x}\C[M]_{y} = M_{x} \otimes 
{}_{x}\C_{y} \overset{\rho_{x,y}}{\rightarrow} M_{y} = \bar{M}_{y}, 
&
\bar{\rho}_{x,M} : \bar{M}_{x} \otimes {}_{x}\C[M]_{M} = M_{x} \otimes 
0 \overset{0}{\rightarrow} k = \bar{M}_{M}, 
\\
\bar{\rho}_{M,x} : \bar{M}_{M} \otimes {}_{M}\C[M]_{x} = k \otimes 
M_{x} \rightarrow M_{x} = \bar{M}_{x}, 
&
\bar{\rho}_{M,M} : \bar{M}_{M} \otimes {}_{M}\C[M]_{M} = k \otimes k 
\rightarrow k = \bar{M}_{M}, 
}
\]
where the last two maps are the action of $k$ on $M_{x}$ and the 
product in $k$, respectively. 
Since $\bar{M} = {}_{M}\C[M]$, we get that $\bar{M}$ is a projective 
$\C[M]$-module satisfying, 
by Yoneda's lemma, $\mathrm{Hom}_{\C[M]} (\bar{M} , \bar{M}) \simeq k$. 
Also, it is easy to see that $\bar{M}$ is small, since 
$\mathrm{Hom}_{\C[M]} (\bar{M} , N) \simeq \mathrm{Hom}_{\C[M]} ({}_{M}\C[M] , N) 
\simeq N_{M}$, for each $\C[M]$-module $N$.

Since $\C$ is a convex subcategory of $\C[M]$ there is a functor $i : 
\mathrm{Mod}_{\C} \rightarrow \mathrm{Mod}_{\C[M]}$, defined at he 
beginning of the this section. 
 
We have that $i({}_{y}\C) = {}_{y}\C[M]$, and hence $i$ preserves 
relatively projectives 
and homotopically projectives, by definition. 

Next we consider the functors $\mathrm{Hom}_{\C[M]} (i(\place),\bar{M}) 
: \mathrm{Mod}_{\C} \rightarrow \mathrm{Mod}_{k}$ and 
$\mathrm{Hom}_{\C} (\place,M) : \mathrm{Mod}_{\C} \rightarrow \mathrm{Mod}_{k}$. 

We remark that they are isomorphic, i.e., there exist a natural 
isomorphism 
\begin{equation}
\label{eq:1}
   \mathrm{Hom}_{\C[M]} (i(\place),\bar{M}) \simeq \mathrm{Hom}_{\C} 
(\place,M),
\end{equation} 
given by 
\begin{align*}
\alpha : \mathrm{Hom}_{\C[M]} (i(\place),\bar{M}) &\rightarrow 
\mathrm{Hom}_{\C} (\place,M)
\\
\{ t_{\bar{x}} \}_{\bar{x} \in \C[M]_{0}} &\mapsto \{ t_{x} \}_{x \in 
\C_{0}},
\end{align*}
with inverse 
\begin{align*}
\beta : \mathrm{Hom}_{\C} (\place,M) &\rightarrow \mathrm{Hom}_{\C[M]} 
(i(\place),\bar{M}) 
\\
\{ t_{x} \}_{x \in \C_{0}} &\mapsto \{ t_{x} \}_{x \in \C_{0}} \sqcup 
\{ 0_{M} \}.
\end{align*}

Since $i$ is an exact functor that preserves injectives, we have that 
$\mathrm{Ext}^{\bullet}_{\C[M]} (i(\place), \bar{M})$ is a universal 
$\delta$-functor, and it is isomorphic in degree zero to 
$\mathrm{Hom}_{\C} (\place, M)$, so 
there is an isomorphism of $\delta$-functors 
\[     \mathrm{Ext}^{\bullet}_{\C[M]} (i(\place), \bar{M}) \simeq 
\mathrm{Ext}^{\bullet}_{\C} (\place, M).     \]

Also, the following identity holds 
\begin{equation}
\label{eq:2}
   \mathrm{Hom}_{\C[M]} (\bar{M} , i(\place) ) = 0.
\end{equation}

\begin{teo}
\label{thm:derone}
Let $\C$ and $\D$ be two $k$-linear categories, $M$ a right $\C$-module 
and $N$ a right $\D$-module. 
For any triangulated equivalence $\phi$ from $D(\D)$ to $D(\C)$, which 
maps $N$ to $M$, 
there exists a triangulated equivalence $\Phi$ from $D(\D[N])$ to 
$D(\C[M])$ which restricts to $\phi$. 
\end{teo}
\noindent\textbf{Proof.} 
According to Theorem \ref{thm:Kel9.2}, $\phi$ is determined by its 
restriction, which is also an equivalence, 
\begin{align*}
\phi' : \D &\rightarrow \E \subset D(\C)
\\
\phi'(y) &= {}_{y}T,
\end{align*}
where ${}_{y}T = \phi({}_{y}\D[0])$ denotes a complex of right 
$\C$-modules ($y \in \D_{0}$). 
By definition of equivalence these complexes form a set of small 
generators of $D(\C)$, such that  
\[     \mathrm{Hom}_{D(\C)} ({}_{y}T , {}_{y'}T[n]) = 0     \]
for $n \neq 0$, and 
\[     \mathrm{Hom}_{D(\C)} ({}_{y}T , {}_{y'}T) = {}_{y'}\D_{y}.     
\]

We are going to define an equivalence from $\D[N]$ to a subcategory 
$\bar{\E}$ of $D(\C[M])$ satisfying the 
hypotheses of Theorem \ref{thm:Kel9.2}. 
The following functor $\Phi'$ is fully faithful 
\begin{align*}
\Phi' : \D[N] &\rightarrow \bar{\E} \subset D(\C[M])
\\
\Phi'(y) &= i({}_{y}T), \hskip 0.5cm \text{if $y \in \D_{0}$,}
\\
\Phi'(N) &= \bar{M}[0].
\end{align*}
The definition on the morphisms is the natural one but it may be useful 
to precise it. 

Let us take $\Phi'(f) = i \circ \phi'(f)$, for $f \in {}_{y'}\D_{y} = 
{}_{y'}\D[N]_{y}$, and $\Phi'(f) = 0$, for 
$f \in {}_{y}\D[N]_{N}$. 
Given $f \in {}_{N}\D[N]_{y}$, we define $\Phi'(f)$ by the following 
chain of natural isomorphisms 
\begin{gather*}
   {}_{N}\D[N]_{y} = N_{y} \overset{\simeq}{\rightarrow} 
\mathrm{Hom}_{\D} ({}_{y}\D , N) 
   \overset{inc}{\rightarrow} \mathrm{Hom}_{D(\D)} ({}_{y}\D , N[0]) 
   \\
   \overset{\phi}{\rightarrow} \mathrm{Hom}_{D(\C)} ({}_{y}T , M[0]) 
   \overset{\beta'}{\rightarrow} \mathrm{Hom}_{D(\C[M])} (i({}_{y}T) , 
\bar{M}[0]) 
   \overset{\simeq}{\rightarrow} \mathrm{Hom}_{\mathcal{H}(\C[M])} 
(i({}_{y}T) , \bar{M}[0]),
\end{gather*}             
where $\beta'$ is the morphism induced by $\beta$ on 
$\mathcal{H}_{p}^{b}$. 
We remark that the last isomorphism holds since $\bar{M}$ is 
$\C[M]$-projective. 
It remains to check that $\beta'$ is an isomorphism: taking into 
account the short exact sequence \eqref{eq:seq}, 
one only needs to check that it is so on each ${}_{y}\C[n]$. 
This is quite simple and follows from the isomorphisms 
\[        \mathrm{Hom}_{D(\C)} ({}_{y}\C[0] , M[0]) = \mathrm{Hom}_{\C} 
({}_{y}\C , M) 
          \overset{\beta}{\rightarrow} \mathrm{Hom}_{\C[M]} 
(i({}_{y}\C) , \bar{M}),     \]
and 
\begin{gather*}
        \mathrm{Hom}_{D(\C)} ({}_{y}\C[n] , M[0]) = 
\mathrm{Hom}_{D(\C)} ({}_{y}\C[0] , M[-n]) 
          \simeq \mathrm{Ext}^{-n}_{\C} ({}_{y}\C , M) = 0 
   \\
   \rightarrow \mathrm{Hom}_{D(\C[M])} (i({}_{y}\C)[0] , \bar{M}[-n]) = 
\mathrm{Hom}_{D(\C[M])} (i({}_{y}\C[n]), \bar{M}[0]),     
\end{gather*}
for $n \neq 0$.
The last map is an isomorphism since $i({}_{y}\C) = {}_{y}\C[M]$ and, 
for $n \neq 0$, we have that 
\[     \mathrm{Hom}_{D(\C[M])} (i({}_{y}\C)[0] , \bar{M}[n]) = 
\mathrm{Ext}^{-n}_{\C[M]} (i({}_{y}\C) , \bar{M}) = 0.     \]

Finally, for $f \in {}_{N}\D[N]_{N}$, we define $\Phi'(f)$ by means of 
the isomorphisms 
\[     {}_{N}\D[N]_{N} = k \simeq \mathrm{Hom}_{\C[M]} (\bar{M} , 
\bar{M}) 
       = \mathrm{Hom}_{D(\C[M])} (\bar{M}[0] , \bar{M}[0]).     \]

The functor $\Phi$ is fully faithful by definition. 
Since $i$ is fully faithful and preserves homotopically projectives, 
\begin{gather*}
   \mathrm{Hom}_{D(\C[M])} (i({}_{y}T) , i({}_{y'}T)) = 
\mathrm{Hom}_{\mathcal{H}(\C[M])} (i({}_{y}T) , i({}_{y'}T)) 
   = \mathrm{Hom}_{\mathcal{H}(\C)} ({}_{y}T , {}_{y'}T) = 
\mathrm{Hom}_{D(\C)} ({}_{y}T , {}_{y'}T) = {}_{y'}\D_{y},    
\end{gather*}
for $y, y'\in \D_{0}$. 
Also, $\mathrm{Hom}_{D(\C[M])} (\bar{M} , i({}_{y'}T)) = 0$ as a 
consequence of \eqref{eq:2}.
All other cases are straightforward. 

We also need to prove that $\mathrm{Hom}_{D(\C[M])} (\Phi'(\bar{y}) , 
\Phi'(\bar{y}')[n]) = 0$, for $n \neq 0$. 
This is achieved in exactly the same way as before, just considering a 
shift by $n$ and noticing that $i$ commutes 
with the shift by definition. 

The image of the functor $\Phi'$ is a set of small generators: 
they are small since $\bar{M}[0]$ is small and $\{ i({}_{y}T) \}_{y \in 
\D_{0}}$ is set of small objects. 
The latter is proved directly from the sequence \eqref{eq:seq} and the 
fact that $i({}_{x}\C) = {}_{x}\C[M]$ is small. 

To prove that they are a set of generators we proceed as follows: 
$\{ {}_{y}T \}_{y \in \D_{0}}$ is a set of generators of $D(\C)$, then 
the full strictly triangulated subcategory 
closed under direct sums containing them also contains $\{ {}_{y}\C 
\}_{y \in \D_{0}}$. 
So, the triangulated subcategory generated by $\{ i({}_{y}T) \}_{y \in 
\D_{0}}$ contains 
$\{ i({}_{y}\C) \}_{y \in \D_{0}}= \{{}_{y}\C[M] \}_{y \in \D_{0}}$. 
As a consequence, the triangulated subcategory generated by the image 
of $\Phi'$ 
contains $\{{}_{y}\C[M] \}_{y \in \D_{0}}$ and $\bar{M} = {}_{M}\C[M]$, 
whence it is the whole $D(\C[M])$. 
\qed

\section{Happel's cohomological long exact sequence}

In this section we first generalize in Theorem \ref{thm:happel} the long exact sequence in 
\cite{Hap1}, Thm. 5.3. to Hochschild-Mitchell cohomology. 
Although the proof is quite similar to the algebraic case but a little 
bit more technical, it is interesting to remark that in the categorical 
context, a more general statement (Thm. \ref{thm:longseq}) not only 
holds but it is in fact more natural. 
The proof of this general statement has been inspired by an article by 
Cibils (cf. Thm. 4.5 \cite{Cib1}) and in fact provides a simpler proof 
to Cibils' result. 

We first state some definitions. 
Given a $\C$-bimodule $N$, let $j(N)$ be the $\C[M]$-bimodule, such 
that 
${}_{\bar{x}}j(N)_{M} = {}_{M}j(N)_{\bar{x}} = 0$, for $\bar{x} \in 
\C[M]_{0}$, and ${}_{y}j(N)_{x} = {}_{y}N_{x}$, for 
$x, y \in \C_{0}$. 
The action is induced by the action of $\C$ on $N$. 
Also, we will denote by $S$ the simple right $\C[M]$-module satisfying 
$S_{x} = 0$, for $x \in \C_{0}$, and $S_{M} = k$. 
The action is the obvious one. 

\begin{lema}
\label{lem:1}
Let $\C$ be a $k$-linear category and $M$ a right $\C$-module. 
The following holds: 
\begin{enumerate}
\item
\label{1} ${}_{}\C[M]^{e} \simeq \C[M]_{M} \otimes_{k} {}_{M}\C[M] 
\simeq \mathrm{Hom}_{k}(S , \bar{M})$, as $\C[M]$-bimodules.  

\item
\label{2} $\mathrm{Ext}^{n+1}_{\C[M]} (S , \bar{M}) \simeq 
\mathrm{Ext}^{n}_{\C} (M , M)$, for $n \geq 1$.

\item
\label{3} $\mathrm{Ext}^{1}_{\C[M]} (S , \bar{M}) \simeq 
\mathrm{Hom}_{\C} (M , M)/k$. 

\item
\label{4} $\mathrm{Hom}_{\C[M]} (S , \bar{M}) = 0$. 

\item
\label{5} $\mathrm{Ext}^{n}_{\C^{e}} (\C , \C) \simeq 
\mathrm{Ext}^{n}_{\C[M]^{e}} (j(\C) , j(\C))$, for $n \geq 0$. 
\end{enumerate} 
\end{lema}
\noindent\textbf{Proof.} 
(\ref{1}). It is clear that the following morphism of $\C[M]$-bimodules
\begin{align*}
{}_{\bar{x}}\phi_{\bar{y}} : \mathrm{Hom}_{k} (S_{\bar{x}} , 
\bar{M}_{\bar{y}}) 
                                         &\rightarrow 
{}_{\bar{x}}\C[M]_{M} \otimes_{k} {}_{M}\C[M]_{\bar{y}}
\\
{}_{\bar{x}}\phi_{\bar{y}} &= 0, \hskip 1.55cm \text{if $\bar{x} \neq 
M$} 
\\
{}_{\bar{x}}\phi_{\bar{y}}(f) &= 1 \otimes f(1), \hskip 0.5cm \text{if 
$\bar{x} = M$}, 
f \in \mathrm{Hom}_{k} (k , \bar{M}_{\bar{y}}). 
\end{align*}
is in fact an isomorphism. 

In order to prove (\ref{2}), (\ref{3}) and (\ref{4}) we proceed as 
follows. 
There is a short exact sequence of right $\C[M]$-modules
\[     0 \rightarrow i(M) \overset{f}{\rightarrow} \bar{M} 
\overset{g}{\rightarrow} S \rightarrow 0.     \]
The morphisms are the obvious ones. 
Applying the functor $\mathrm{Hom}_{\C[M]} (\place , \bar{M})$ to this 
short exact sequence 
we get the long exact sequence 
\begin{align*}
      0 &\rightarrow \mathrm{Hom}_{\C[M]}(S , \bar{M}) \rightarrow 
\mathrm{Hom}_{\C[M]} (\bar{M} , \bar{M}) 
         \rightarrow \mathrm{Hom}_{\C[M]} (i(M) , \bar{M}) \rightarrow 
\mathrm{Ext}^{1}_{\C[M]} (S , \bar{M}) \rightarrow \dots
\\
     \dots &\rightarrow \mathrm{Ext}^{n}_{\C[M]} (S , \bar{M}) 
\rightarrow \mathrm{Ext}^{n}_{\C[M]} (\bar{M} , \bar{M}) \rightarrow 
\mathrm{Ext}^{n}_{\C[M]} (i(M) , \bar{M}) \rightarrow \mathrm{Ext}^{n+1}_{\C[M]} 
(S , \bar{M}) \rightarrow \dots
\end{align*}

Taking into account that $i$ preserves exactness and projectives, and 
the isomorphism \eqref{eq:1}, we have that 
$\mathrm{Hom}_{\C[M]} (i(M) , \bar{M}) \simeq \mathrm{Hom}_{\C} (M , 
M)$ and 
$\mathrm{Ext}^{n}_{\C[M]} (i(M) , \bar{M}) \simeq \mathrm{Ext}^{n}_{\C} 
(M , M)$, for $n \geq 1$. 
Also, we see immediately that $\mathrm{Ext}^{n}_{\C[M]} (\bar{M} , 
\bar{M}) = 0$, for $n \geq 1$, 
since $\bar{M} = {}_{M}\C[M]$ is projective. 
This proves (\ref{2}). 

For the other statements, we recall that $\mathrm{Hom}_{\C[M]} (\bar{M} 
, \bar{M}) \simeq {}_{M}\C[M]_{M} = k$, 
and notice that the map given by 
$f^{*} : \mathrm{Hom}_{\C[M]} (\bar{M} , \bar{M}) \rightarrow 
\mathrm{Hom}_{\C[M]} (i(M) , \bar{M})$ is not zero 
since $f^{*} (\id_{\bar{M}}) = f \neq 0$, and so $f^{*}$ is injective. 
Hence we get (\ref{3}) and (\ref{4}). 

In order to prove (\ref{5}) we only use that $\C^{e}$ is a full convex 
subcategory of $\C[M]^{e}$ and apply Lemma \ref{lema:convex}. 
\qed

\begin{teo}
\label{thm:happel}
Let $\C$ be a $k$-linear category and $M$ a right $\C$-module. 
There is a cohomological long exact sequence 
\begin{gather*}
      0 \rightarrow HH^{0} (\C[M]) \rightarrow HH^{0} (\C) \rightarrow 
\mathrm{Hom}_{\C} (M , M)/k 
          \rightarrow HH^{1} (\C[M]) \rightarrow HH^{1} (\C) 
\rightarrow \mathrm{Ext}^{1}_{\C} (M , M)
          \rightarrow \dots
\\
          \dots \rightarrow \mathrm{Ext}^{n-1}_{\C} (M , M) \rightarrow 
HH^{n} (\C[M]) \rightarrow HH^{n} (\C) \rightarrow 
\mathrm{Ext}^{n}_{\C} (M , M) \rightarrow HH^{n+1} (\C[M])\rightarrow \dots
\end{gather*}
\end{teo}
\noindent\textbf{Proof.}
Let us consider the following short exact sequence of $\C[M]$-bimodules 
\begin{equation}
\label{eq:seq2}
     0 \rightarrow K \overset{\alpha}{\rightarrow} \C[M] 
\overset{\beta}{\rightarrow} j(\C) \rightarrow 0,     
\end{equation}
where $\beta$ is given by $\beta(f) = f$, for $f \in {}_{y}\C_{x} 
\subset {}_{y}\C[M]_{x}$, and 
zero in any other case. 
The $\C[M]$-bimodule $K$ is its kernel. 

We shall see that $K$ and $\C[M]_{M} \otimes_{k} {}_{M}\C[M]$ are 
isomorphic as $\C[M]$-bimodules. 
To prove this fact we proceed as follows: consider the map 
\begin{align*}
\gamma : \C[M]_{M} \otimes_{k} {}_{M}\C[M] &\rightarrow \C[M]
\\
\gamma(c \otimes c') &= c.c'.
\end{align*} 
It is evident that $\beta \circ \gamma = 0$ and that $\gamma$ is a 
monomorphism. 
If $c \in \mathrm{Ker}(\beta)$, then either $c = 0$ or $c \in 
{}_{M}\C[M]_{\bar{x}}$. 
In this case, $c = \gamma({}_{M}1_{M} \otimes c)$, and hence $c \in 
\mathrm{Im}(\gamma)$. 
It follows that $\C[M]_{M} \otimes_{k} {}_{M}\C[M]$ is also a kernel of 
$\beta$. 
As a consequence, $\mathrm{Ext}^{n}_{\C[M]^{e}} (K , j(\C)) = 0$ for $n 
\geq 1$. 
Also $\mathrm{Hom}_{\C[M]^{e}} (K , j(\C)) = {}_{M}j(\C)_{M} = 0$, so 
$\mathrm{Ext}^{n}_{\C[M]^{e}} (K , j(\C)) = 0$, 
for $n \geq 0$. 

Now, applying the functor $\mathrm{Hom}_{\C[M]^{e}} (\place , j(\C))$ 
to the sequence \eqref{eq:seq2} 
and using that $\mathrm{Ext}^{n}_{\C[M]^{e}} (K , j(\C)) = 0$ for $n 
\geq 0$, we get 
$\mathrm{Ext}^{n}_{\C[M]^{e}} (j(\C) , j(\C)) = 
\mathrm{Ext}^{n}_{\C[M]^{e}} (\C[M] , j(\C))$, 
for $n \geq 0$. 
The first one is isomorphic to $HH^{n}(\C)$ using Lemma \ref{lem:1}, 
(\ref{5}). 

Also, 
\[     H^{n} (\C[M] , K) = \mathrm{Ext}^{n}_{\C[M]^{e}} (\C[M] , K) 
\simeq \mathrm{Ext}^{n}_{\C[M]} (S , \bar{M}) 
       = H^{n} (\C[M] , \mathrm{Hom}_{k} (S , \bar{M})),     \]
by Lemma \ref{lem:1}, (\ref{1}). 
Finally, $H^{n}(\C[M], \mathrm{Hom}_{k}(S,\bar{M}))$ is isomorphic to 
$\mathrm{Ext}_{\C[M]}^{n} (S , \bar{M})$ since, 
by adjunction, the complex computing the Hochschild-Mitchell cohomology 
also gives the $\mathrm{Ext}$ groups. 
We also notice that $\mathrm{Hom}_{\C[M]^{e}} (\C[M] , K) = 0$, 
$\mathrm{Ext}_{\C[M]^{e}}^{1} (\C[M] , K) = \mathrm{Hom}_{\C}(M , M)/k$ 
and 
$\mathrm{Ext}_{\C[M]^{e}}^{n} (\C[M] , K) = \mathrm{Ext}^{n-1}_{\C}(M , 
M)$, for $n \geq 2$, using 
Lemma \ref{lem:1}, (\ref{4}), (\ref{3}) and (\ref{1}) respectively. 

Applying now the functor $\mathrm{Hom}_{\C[M]^{e}} (\C[M] , \place)$ to 
\eqref{eq:seq2}, we obtain the long exact sequence 
\begin{align*}
     0 &\rightarrow \mathrm{Hom}_{\C[M]^{e}} (\C[M] , K) \rightarrow 
\mathrm{Hom}_{\C[M]^{e}} (\C[M] , \C[M]) 
        \rightarrow \mathrm{Hom}_{\C[M]^{e}} (\C[M] , j(\C)) 
\rightarrow \mathrm{Ext}^{1}_{\C[M]^{e}} (\C[M] , K) 
        \rightarrow \dots
\\
        \dots &\rightarrow \mathrm{Ext}^{n}_{\C[M]^{e}} (\C[M] , K) 
\rightarrow \mathrm{Ext}^{n}_{\C[M]^{e}} (\C[M] , \C[M]) 
            \rightarrow \mathrm{Ext}^{n}_{\C[M]^{e}} (\C[M] , j(\C)) 
\rightarrow \mathrm{Ext}^{n+1}_{\C[M]^{e}} (\C[M] , K) 
            \rightarrow \dots
\end{align*}
Using the identifications above the theorem follows.
\qed

Next we will consider a more general situation. 
Let $\C_{1}$ and $\C_{2}$ be two $k$-linear categories, and let $M$ be 
a $\C_{1}$-$\C_{2}$-bimodule. 
We define the category $\C = \C_{1} \sqcup_{M} \C_{2}$ with objects 
$\C_{0} = (\C_{1})_{0} \sqcup (\C_{2})_{0}$ 
and morphisms 
\[
{}_{x}\C_{y} = \begin{cases} 
                       {}_{x}(\C_{1})_{y}, \hskip 0.5cm \text{for $x,y 
\in (\C_{1})_{0}$,} 
                       \\
                       {}_{x}(\C_{2})_{y}, \hskip 0.5cm \text{for $x,y 
\in (\C_{2})_{0}$,} 
                       \\
                       {}_{x}M_{y}, \hskip 0.5cm \text{for $x \in 
(\C_{1})_{0}$, $y \in (\C_{2})_{0}$,} 
                       \\
                       0, \hskip 0.5cm \text{otherwise.}                        
                       \end{cases}
\]
\begin{ejem}
If $(\C_{1})_{0} = \{*\}$, ${}_{*}(\C_{1})_{*} = k$ and $M$ is a right 
$\C_{2}$-module, then 
$\C_{1} \sqcup_{M} \C_{2} = \C_{2}[M]$. 
\end{ejem}

Since, for $i,j \in \{ 1,2 \}$, $\C_{i} \boxtimes \C_{j}^{op}$ is a 
convex subcategory of $\C^{e}$, there are well-defined restriction functors. 
Given a $\C$-bimodule $N$, we shall denote $r_{i,j}(N)$ the 
corresponding restriction. 
We also write $r_{i}(N) = r_{i,i}(N)$. 

In this situation there is a cohomological long exact sequence 
generalizing the previous one. 
The key fact of the proof is that it is possible to decompose the 
Hochschild-Mitchell projective resolution of $\C$ as 
$\C$-bimodule as follows  
\begin{align*}
   N_{n}(\C) 
   &= \bigoplus_{\hbox{\footnotesize $(x_{0}, \dots, x_{n}) \in 
\C_{0}^{n+1}$}} {}_{\place}\C_{x_{n}} \otimes {}_{x_{n}}\C_{x_{n-1}} \otimes 
\dots  
   \otimes {}_{x_{1}}\C_{x_{0}} \otimes {}_{x_{0}}\C_{\place} 
   \\
   &= \bigoplus_{\hbox{\footnotesize $(x_{0}, \dots, x_{n}) \in 
(\C_{1})_{0}^{n+1}$}} {}_{\place}\C_{x_{n}} \otimes 
{}_{x_{n}}(\C_{1})_{x_{n-1}} \otimes \dots \otimes {}_{x_{1}}(\C_{1})_{x_{0}} \otimes 
{}_{x_{0}}\C_{\place} 
   \\
   &\oplus 
     \bigoplus_{\hbox{\footnotesize $(x_{0}, \dots, x_{n}) \in 
(\C_{2})_{0}^{n+1}$}} {}_{\place}\C_{x_{n}} \otimes 
{}_{x_{n}}(\C_{2})_{x_{n-1}} \otimes \dots \otimes {}_{x_{1}}(\C_{2})_{x_{0}} \otimes 
{}_{x_{0}}\C_{\place}    
   \\
    &\oplus 
    \bigoplus_{i = 0}^{n-1}\bigoplus_{\hbox{\footnotesize 
$\begin{matrix}(x_{0}, \dots, x_{i}) \in (\C_{2})_{0}^{i+1}\\ 
                                                      (x_{i+1}, \dots, 
x_{n}) \in (\C_{1})_{0}^{n-i}\end{matrix}$}} 
          {}_{\place}\C_{x_{n}} \otimes {}_{x_{n}}(\C_{1})_{x_{n-1}} 
\otimes \dots \otimes {}_{x_{i+1}}M_{x_{i}} 
          \otimes \dots \otimes {}_{x_{1}}(\C_{2})_{x_{0}} \otimes 
{}_{x_{0}}\C_{\place}.    
\end{align*}
This decomposition gives 
\begin{align*}
   \mathrm{Hom}_{\C^{e}} (N_{n}(\C) , N) 
   &= \prod_{\hbox{\footnotesize $(x_{0}, \dots, x_{n}) \in 
\C_{0}^{n+1}$}} 
        \mathrm{Hom}_{k} ({}_{x_{n}}\C_{x_{n-1}} \otimes \dots \otimes 
{}_{x_{1}}\C_{x_{0}} , {}_{x_{n}}N_{x_{0}}) 
   \\
   &= \prod_{\hbox{\footnotesize $(x_{0}, \dots, x_{n}) \in 
(\C_{1})_{0}^{n+1}$}} 
   \mathrm{Hom}_{k} ({}_{x_{n}}(\C_{1})_{x_{n-1}} \otimes \dots \otimes 
{}_{x_{1}}(\C_{1})_{x_{0}} , {}_{x_{n}}N_{x_{0}}) 
   \\
   &\oplus 
   \prod_{\hbox{\footnotesize $(x_{0}, \dots, x_{n}) \in 
(\C_{2})_{0}^{n+1}$}} 
   \mathrm{Hom}_{k} ({}_{x_{n}}(\C_{2})_{x_{n-1}} \otimes \dots \otimes 
{}_{x_{1}}(\C_{2})_{x_{0}} , {}_{x_{n}}N_{x_{0}})    
   \\
   &\oplus 
   \bigoplus_{i = 0}^{n-1}\prod_{\hbox{\footnotesize 
$\begin{matrix}(x_{0}, \dots, x_{i}) \in (\C_{2})_{0}^{i+1}\\ 
                                                      (x_{i+1}, \dots, 
x_{n}) \in (\C_{1})_{0}^{n-i}\end{matrix}$}} 
          \mathrm{Hom}_{k} ({}_{x_{n}}(\C_{1})_{x_{n-1}} \otimes \dots 
\otimes {}_{x_{i+1}}M_{x_{i}} 
          \otimes \dots \otimes {}_{x_{1}}(\C_{2})_{x_{0}} , 
{}_{x_{n}}N_{x_{0}}) 
   \\
   &= \mathrm{Hom}_{\C_{1}^{e}} (N_{n}(\C_{1}) , r_{1}(N)) \oplus 
      \mathrm{Hom}_{\C_{2}^{e}} (N_{n}(\C_{2}) , r_{2}(N)) 
      \oplus
      \mathrm{Hom}_{\C_{1} \boxtimes \C_{2}^{op}} (\tilde{M}_{n-1} , 
r_{1,2}(N)),    
\end{align*}
where ($\tilde{M}_{n}$, $d_{n}$) is the complex of projective 
$\C_{1}$-$\C_{2}$-bimodules given by 
\[         \tilde{M}_{n} = 
          \bigoplus_{i = 0}^{n}\bigoplus_{\hbox{\footnotesize 
$\begin{matrix}(x_{0}, \dots, x_{i}) \in (\C_{1})_{0}^{i+1}\\ 
                                                      (x_{i+1}, \dots, 
x_{n+1}) \in (\C_{1})_{0}^{n+1-i}\end{matrix}$}} 
          {}_{\place}(\C_{1})_{x_{n+1}} \otimes 
{}_{x_{n+1}}(\C_{2})_{x_{n}} \otimes \dots \otimes {}_{x_{i+1}}M_{x_{i}} 
          \otimes \dots \otimes {}_{x_{1}}(\C_{2})_{x_{0}} \otimes 
{}_{x_{0}}(\C_{2})_{\place},    
\]
with differential $d_{\bullet}$ obtained by restricting the 
differential of the Hochschild-Mitchell resolution. 
This complex is in fact a projective resolution of $M$ as a 
$\C_{1}$-$\C_{2}$-bimodule. 
In order to prove this statement, it is sufficient to notice that 
($\tilde{M}_{\bullet}$, $d_{\bullet}$) 
is the total complex obtained from the first quadrant double complex 
\[     \tilde{M}_{i,j} = \bigoplus_{\hbox{\footnotesize 
$\begin{matrix}(x_{0}, \dots, x_{i}) \in (\C_{2})_{0}^{i+1}\\ 
                                                      (x_{i+1}, \dots, 
x_{i+j+1}) \in (\C_{1})_{0}^{j}\end{matrix}$}} 
          {}_{\place}(\C_{1})_{x_{n+1}} \otimes 
{}_{x_{n+1}}(\C_{1})_{x_{n}} \otimes \dots \otimes {}_{x_{i+1}}M_{x_{i}} 
          \otimes \dots \otimes {}_{x_{1}}(\C_{2})_{x_{0}} \otimes 
{}_{x_{0}}(\C_{2})_{\place},     \]
where the vertical and horizontal differentials are 
\begin{align*}
{}_{y}(d^{h}_{i,j})_{x} 
&({}_{y}(c_{1})_{x_{n+1}} \otimes \dots \otimes {}_{x_{i+1}}m_{x_{i}} 
\otimes \dots \otimes {}_{x_{0}}(c_{2})_{x}) 
\\
= &{}_{y}(c_{1})_{x_{n+1}}.{}_{x_{n+1}}(c_{1})_{x_{n}} \otimes \dots 
\otimes {}_{x_{i+1}}m_{x_{i}} \otimes \dots \otimes {}_{x_{0}}(c_{2})_{x} 
\\
+ &\sum_{j = i+2}^{n+1} (-1)^{j+n+1} 
{}_{y}(c_{1})_{x_{n+1}} \otimes \dots \otimes 
{}_{x_{j+1}}(c_{1})_{x_{j}}.{}_{x_{j}}(c_{1})_{x_{j-1}} \otimes 
\dots \otimes {}_{x_{i+1}}m_{x_{i}} \otimes \dots \otimes 
{}_{x_{0}}(c_{2})_{x} 
\\
+ &(-1)^{i+n}
{}_{y}(c_{1})_{x_{n+1}} \otimes \dots \otimes 
{}_{x_{i+2}}(c_{1})_{x_{i+1}}.{}_{x_{i+1}}m_{x_{i}} \otimes \dots \otimes {}_{x_{0}}(c_{2})_{x} 
\end{align*}
and 
\begin{align*}
{}_{y}(d^{v}_{i,j})_{x} 
&({}_{y}(c_{1})_{x_{n+1}} \otimes \dots \otimes {}_{x_{i+1}}m_{x_{i}} 
\otimes \dots \otimes {}_{x_{0}}(c_{2})_{x}) 
\\
= &(-1)^{i+n+1}
{}_{y}(c_{1})_{x_{n+1}} \otimes \dots \otimes 
{}_{x_{i+1}}m_{x_{i}}.{}_{x_{i}}(c_{2})_{x_{i-1}}  \otimes \dots \otimes {}_{x_{0}}(c_{2})_{x} 
\\
+ &\sum_{j = 1}^{i-1} (-1)^{j+n+1}
{}_{y}(c_{1})_{x_{n+1}} \otimes \dots \otimes {}_{x_{i+1}}m_{x_{i}} 
\otimes 
\dots \otimes {}_{x_{j+1}}(c_{2})_{x_{j}}.{}_{x_{j}}(c_{2})_{x_{j-1}} 
\otimes \dots \otimes {}_{x_{0}}(c_{2})_{x} 
\\
+ &(-1)^{n+1}
{}_{y}(c_{1})_{x_{n+1}} \otimes \dots \otimes {}_{x_{i+1}}m_{x_{i}} 
\otimes \dots \otimes {}_{x_{1}}(c_{2})_{x_{0}}.{}_{x_{0}}(c_{2})_{x}. 
\end{align*}
This double complex has exact rows and columns using the usual homotopy 
arguments. 
Then the cohomology of the cochain complex 
($\mathrm{Hom}_{\C_{1} \boxtimes \C_{2}^{op}} (\tilde{M}_{\bullet} , 
r_{1,2}(N))$, $d^{*}_{\bullet}$) is exactly 
$\mathrm{Ext}^{\bullet}_{\C_{1} \boxtimes \C_{2}^{op}} (M , 
r_{1,2}(N))$. 

We also notice that this cochain complex is actually a subcomplex of 
$\mathrm{Hom}_{\C^{e}} (N_{n}(\C) , N)$, which is the complex computing the 
Hochschild-Mitchell cohomology of $\C$, and its quotient is 
$\mathrm{Hom}_{\C_{1}^{e}} (N_{n}(\C_{1}) , r_{1}(N)) \oplus 
\mathrm{Hom}_{\C_{2}^{e}} (N_{n}(\C_{2}) , r_{2}(N))$. 
In other words, there is a short exact sequence of complexes of 
$k$-modules 
\begin{align*}
     0 &\rightarrow 
       \mathrm{Hom}_{\C_{1} \boxtimes \C_{2}^{op}} 
(\tilde{M}_{\bullet-1} , r_{1,2}(N)) 
       \rightarrow \mathrm{Hom}_{\C^{e}} (N_{\bullet}(\C) , N) 
\\
&\rightarrow 
       \mathrm{Hom}_{\C_{1}^{e}} (N_{\bullet}(\C_{1}) , r_{1}(N)) 
       \oplus \mathrm{Hom}_{\C_{2}^{e}} (N_{\bullet}(\C_{2}) , 
r_{2}(N)) 
       \rightarrow 0. 
\end{align*} 

The cohomological long exact sequence obtained from this short exact 
sequence yields the following theorem
\begin{teo}
\label{thm:longseq}
Let $\C_{1}$ and $\C_{2}$ be two small $k$-linear categories, and let 
$M$ be a $\C_{1}$-$\C_{2}$-bimodule. 
Denoting $\C = \C_{1} \sqcup_{M} \C_{2}$, there is cohomological long 
exact sequence  
\begin{align*}
     0 &\rightarrow H^{0} (\C , N) \rightarrow H^{0}(\C_{1} , r_{1}(N)) 
\oplus H^{0} (\C_{2} , r_{2}(N)) 
           \rightarrow \mathrm{Hom}_{\C_{1} \boxtimes \C_{2}^{op}} (M , 
r_{1,2}(N)) \rightarrow H^{1} (\C , N) \rightarrow \dots 
\\
        \dots &\rightarrow H^{n} (\C , N) \rightarrow H^{n}(\C_{1} , 
r_{1}(N)) \oplus H^{n} (\C_{2} , r_{2}(N)) 
           \rightarrow \mathrm{Ext}_{\C_{1} \boxtimes \C_{2}^{op}}^{n} 
(M , r_{1,2}(N)) \rightarrow H^{n+1} (\C , N) \rightarrow \dots 
\end{align*}
\end{teo}
This theorem provides a long exact sequence generalizing the one 
obtained by Cibils \cite{Cib1} for algebras and the one point-extension 
sequence proved before. 



\footnotesize 
\noindent E.H.:
\\Departamento de Matem\'atica,
 Facultad de Ciencias Exactas y Naturales,\\
 Universidad de Buenos Aires
\\Ciudad Universitaria, Pabell\'on 1\\
1428, Buenos Aires, Argentina. \\{\tt eherscov@bigua.dm.uba.ar}

\medskip 

\noindent A.S.:
\\Departamento de Matem\'atica,
 Facultad de Ciencias Exactas y Naturales,\\
 Universidad de Buenos Aires
\\Ciudad Universitaria, Pabell\'on 1\\
1428, Buenos Aires, Argentina. \\{\tt asolotar@dm.uba.ar}

%

\end{document}